\documentclass[10 pt,reqno]{amsart}
\usepackage{amsfonts}
\usepackage{amsmath}
\usepackage{amsthm}
\usepackage{amssymb}
\newcommand{\C}{\mathbb{C}}
\newcommand{\Cn}{\mathbb{C}^n}
\newcommand{\R}{\mathbb{R}}
\newcommand{\M}{\mathcal{M}}

\newcommand{\Ran}{\operatorname{Ran}}

\newcommand{\rank}{\operatorname{Rank}}
\newcommand{\Span}{\operatorname{span}}
\newcommand{\ip}[1]{\langle #1 \rangle}
\newcommand{\norm}[1]{\lVert #1 \rVert}
\newcommand{\abs}[1]{\lvert #1 \rvert}

\newtheorem{lemma}{Lemma}
\newtheorem{prop}{Proposition}
\newtheorem{cor}{Corollary}
\newtheorem{theorem}{Theorem}
\newtheorem{conjecture}{Conjecture}
\theoremstyle{definition}
\newtheorem*{defn}{Definition}
\newtheorem{ex}{Example}
\newtheorem*{alg}{Algorithm}

\theoremstyle{remark}
\newtheorem*{rmk}{Remark}
\begin{document}
\title[Unitary equivalence to a complex symmetric matrix]{Unitary equivalence to a complex symmetric matrix: an algorithm}
    \author{James E. Tener}
    \address{Department of Mathematics\\
            Pomona College '08\\
            Claremont, California\\
            91711 \\ USA}
    \email{james.tener@pomona.edu}


    \begin{abstract}
       We present a necessary and sufficient condition for a $3 \times 3$ matrix to be unitarily equivalent to a symmetric matrix with complex entries, and an algorithm whereby an arbitrary $3 \times 3$ matrix can be tested.  This test generalizes to a necessary and sufficient condition that applies to almost every $n \times n$ matrix.  The test is constructive in that it explicitly exhibits the unitary equivalence to a complex symmetric matrix.
    \end{abstract}

\maketitle

\section{Introduction}
\begin{defn}
A \emph{complex symmetric matrix} is an $n \times n$ matrix $T$ with complex entries such that $T = T^t$.  We will often refer to such a matrix as a \emph{CSM}.  We will often refer to a matrix that is unitarily equivalent to a complex symmetric matrix as being \emph{UECSM}.
\end{defn}

\begin{defn}
A \emph{conjugation} is a conjugate-linear operator $C : \Cn \to \Cn$ which is both \emph{involutive} ($C^2 = I$) and \emph{isometric} ($\ip{x,y} = \ip{Cy,Cx}$ for all $x,y \in \Cn$).  We say that a linear operator $T$ on $\Cn$ is $C$-symmetric if $T = CT^*C$.
\end{defn}

The matrices that are unitarily equivalent to complex symmetric matrices can be characterized in terms of $C$-symmetry.  We recall a consequence of \cite[Prop. 2]{GarciaPutinar}:
\begin{prop}\label{Csymm}
An  $n \times n$ matrix $T$ is unitarily equivalent to a complex symmetric matrix if and only if there exists a conjugation $C$ for which the linear operator induced by $T$ is $C$-symmetric.
\end{prop}

For more on complex symmetric operators, see \cite{GarciaPutinar, GarciaPutinar2}.  It is already known that all $2 \times 2$ matrices are unitarily equivalent to a CSM (\cite[Ex. 6]{GarciaPutinar}, \cite[Cor 3.3]{CFT}) but that not all $3 \times 3$ matrices are \cite[Ex. 7]{GarciaPutinar}.  It is also known that all matrices are \emph{similar} to a CSM (\cite[Ex. 4]{GarciaPutinar}, \cite[Thm. 4.4.9]{HornJohnson}), which makes it difficult to determine if a particular matrix is unitarily equivalent to a complex symmetric matrix.  For instance, exactly one of the following matrices is unitarily equivalent to a CSM:
\[
T_1 = \left[\begin{array}{ccc}0 & 7 & 0 \\ 0 & 1 & -5 \\ 0 & 0 & 6\end{array}\right], \quad
T_2 = \left[\begin{array}{ccc}0 & 7 & 0 \\ 0 & 1 & -5 \\ 0 & 0 & 3\end{array}\right].
\]
For more, see Example \ref{intex}.\\

This note presents necessary and sufficient conditions for a $3 \times 3$ matrix to be UECSM, as well as an algorithm with which one can apply these conditions to an arbitrary $3 \times 3$ matrix.  This test also applies to almost every $n \times n$ matrix $T$, and if $T$ is UECSM then it possible to explicitly construct a conjugation $C$ for which $T$ is $C$-symmetric.

\section{Eigenvectors of the Cartesian decomposition}

Any matrix $T$ has a unique Cartesian decomposition $T=A+iB$, where $A$ and $B$ are self-adjoint.  $A$ and $B$ are given explicitly by the formulas $A = \frac12 (T + T^*)$ and $B = \frac1{2i} (T - T^*)$.  It follows that $T$ is $C$-symmetric for some conjugation $C$ if and only if $A$ and $B$ are $C$-symmetric.

We recall a result from \cite[Lemma 1]{GarciaPutinar}:

\begin{lemma}\label{Crealbasis}
If $C$ is a conjugation on $\Cn$, then there is an orthonormal basis $\{e_i\}_{i=1}^n$ such that $Ce_i = e_i$ for all $i$.
\end{lemma}

Such an orthonormal basis is called \emph{$C$-real}.  The following lemma shows that one way such bases arise is from orthonormal bases of eigenvectors of self-adjoint matrices.

\begin{lemma}\label{onbeigenvectors}
If $A$ is a $C$-symmetric self-adjoint matrix, then there exists a $C$-real orthonormal basis of eigenvectors of $A$.
\end{lemma}

\begin{proof}
If $\lambda$ is an eigenvalue of $A$ with eigenvector $x$, then \mbox{$ACx = CAx = C\lambda x = \lambda Cx$}, since $\lambda$ is real.  Since $C$ is involutive, the above equation yields that $C$ must map the eigenspaces of $A$ onto themselves.  Thus we can decompose $C = C_1 \oplus \cdots \oplus C_m$ where the $C_i$ are conjugations on the eigenspaces of $A$ corresponding to distinct eigenvalues.  Since the $i$th eigenspace must have a $C_i$-real orthonormal basis, and the eigenspaces of $A$ are mutually orthogonal, we are done by Lemma \ref{Crealbasis}.
\end{proof}

The following lemma characterizes all UECSM matrices in terms of such orthonormal bases of eigenvectors.

\begin{lemma}\label{realangles}
If $T = A + iB$ is an $n \times n$ matrix, then $T$ is unitarily equivalent to a symmetric matrix if and only if there exist orthonormal bases of eigenvectors $\{e_i\}$ and $\{f_i\}$ of $A$ and $B$ respectively such that $\ip{e_i, f_j} \in \R$ for all $1 \le i,j \le n$.
\end{lemma}
\begin{proof}
First assume that $T$ is UECSM.  By Proposition \ref{Csymm}, there must exist a conjugation $C$ on $\Cn$ such that $CT^*C = T$, and thus $CAC = A$ and $CBC = B$.  By Lemma \ref{onbeigenvectors}, $A$ and $B$ have $C$-real orthonormal bases of eigenvectors $\{e_i\}$ and $\{f_i\}$, respectively.  For these we have,
\[
\ip{e_i, f_j} =  \ip{Cf_j, Ce_i}
 = \ip{f_j, e_i} 
 = \overline{\ip{e_i, f_j}},
\]
and so $\ip{e_i,f_j} \in \R$.\\

Conversely suppose there exist such $\{e_i\}$ and $\{f_i\}$.  Then define $C : \Cn \to \Cn$ by 
\[
Cx = \sum_{i=1}^n \overline{\ip{x, e_i}}e_i.
\]  
It is easy to check that $C$ is a conjugation, and we can calculate
\[
Ce_j = \sum_{i=1}^n \overline{\ip{e_j, e_i}}e_i
 = \sum_{i=1}^n \delta_{i,j}e_i
= e_j
\]
and,
\[
Cf_j = \sum_{i=1}^n \overline{\ip{f_j, e_i}}Ce_i
= \sum_{i=1}^n \ip{f_j, e_i}e_i
= f_j.
\]
Since we know $Ae_i = \lambda_i e_i$ for some $\lambda_i \in \R$, we can calculate,
\[
CACe_i = CAe_i
= \lambda_i e_i
= Ae_i.
\]
By linearity, it follows that $Ax = CACx$ for all $x \in \Cn$.  Similarly one can show that $CBC = B$, and so $T$ is $C$-symmetric.
\end{proof}

\begin{rmk}
In the previous proof, we could have let $C$ be complex conjugation with respect to $\{f_i\}$ instead.
\end{rmk}

Although the following corollary is well-known \cite[Ex. 2.8]{GarciaCCO}, the above lemma provides another proof.
\begin{cor}\label{normal}
Every normal matrix is unitarily equivalent to a complex symmetric matrix.
\end{cor}
\begin{proof}
$T=A+iB$ being normal implies that $A$ and $B$ share a basis of eigenvectors, and so in the hypothesis of Lemma \ref{realangles} we may pick $\{e_i\}$ and $\{f_i\}$ to be the same orthonormal basis.
\end{proof}

\begin{cor}
Let $T = A + iB$ be an $n \times n$ matrix that is unitarily equivalent to a complex symmetric matrix, and let $\{e_i\}$ and $\{f_i\}$ be a pair of orthonormal bases of $A$ and $B$, respectively, as in Lemma \ref{realangles}.  Then for any $A^\prime$ and $B^\prime$ such that the $\{e_i\}$ are eigenvectors for one and $\{f_i\}$ are eigenvectors for the other, \mbox{$T^\prime = A^\prime + iB^\prime$} is unitarily equivalent to a complex symmetric matrix.
\end{cor}

\section{A test for $n$ $\times$ $n$ matrices}

\begin{defn}
We say that a pair of orthogonal bases $\{g_i\}$ and $\{h_i\}$ for $\Cn$ are \emph{proper} if $\ip{g_1,h_1} \in \R$ and $\ip{g_i, h_j} = 0 \implies i\ne1$ and $j\ne1$.  If we let $M=(\ip{g_i,h_j})_{i,j=1}^n$, then this is equivalent to the top-left entry of $M$ being real while the first row and column contain no zeros.
\end{defn}

\begin{theorem}\label{test}
Let $T = A + iB$ be an $n \times n$ matrix, and let $\{g_i\}$ and $\{h_i\}$ be any proper pair of orthogonal bases of eigenvectors of $A$ and $B$, respectively.  Then $T$ is unitarily equivalent to a complex symmetric matrix if for $2 \le i,j \le n$,
\begin{equation}\label{realratio}
\frac{\ip{g_i,h_j}}{\ip{g_i,h_1}\ip{g_1,h_j}} \in \R.
\end{equation}
Moreover, if $A$ and $B$ both have $n$ distinct eigenvalues, then \eqref{realratio} is also necessary for $T$ to be UECSM.
\end{theorem}
\begin{proof}
Suppose that \eqref{realratio} holds.  Define $e_1 = g_1$, $f_1 = h_1$ and otherwise,
\[\begin{array}{cc}
e_i = \frac{1}{\ip{g_i,h_1}}\mbox{ }g_i,& f_j = \frac{1}{\ip{h_j,g_1}}\mbox{ }h_j.
\end{array}\]
Once normalized, these bases satisfy Lemma \ref{realangles}, showing that $T$ is UECSM.\\

Conversely, suppose that $T$ is UECSM and that $A$ and $B$ have $n$ distinct eigenvalues.  By Lemma \ref{realangles}, there are orthonormal bases of eigenvectors $\{e_i\}$ and $\{f_i\}$ of $A$ and $B$, respectively, such that $\ip{e_i, f_j} \in \R$.  As the eigenspaces of $A$ and $B$ are one-dimensional, we can reorder these bases so that $g_i = \omega_ie_i$ and $h_j = \zeta_jf_j$ for unimodular $\omega_i,\zeta_j \in \C$.  Then for $2 \le i,j \le n$,
\begin{eqnarray*}
\frac{\ip{g_i,h_j}\ip{g_1,h_1}}{\ip{g_i,h_1}\ip{g_1,h_j}} & = & \frac{\ip{\omega_ie_i,\zeta_jf_j}\ip{\omega_1e_1,\zeta_1f_1}}{\ip{\omega_ie_i,\zeta_1f_1}\ip{\omega_1e_1,\zeta_jf_j}}\\
& = & \frac{\omega_i\omega_1\overline{\zeta_j\zeta_1}\ip{e_i,f_j}\ip{e_1,f_1}}{\omega_i\omega_1\overline{\zeta_j\zeta_1}\ip{e_i,f_1}\ip{e_1,f_j}}\\
& = & \frac{\ip{e_i,f_j}\ip{e_1,f_1}}{\ip{e_i,f_1}\ip{e_1,f_j}}
\end{eqnarray*}
which is real since each $\ip{e_i, f_j} \in \R$.  The pair $\{g_i\}$ and $\{h_i\}$ being proper ensures that $\ip{g_1, h_1} \in \R \setminus \{0\}$, and so the preceding equation yields that $\{g_i\}$ and $\{h_i\}$ satisfy \eqref{realratio}.
\end{proof}

\begin{rmk}
The condition \eqref{realratio} of Theorem \ref{test} can be visualized using matrices.  If $\{e_i\}$ and $\{f_i\}$ are a proper pair of orthogonal bases of eigenvectors, we can consider the matrix $M = (m_{i,j}) = (\ip{e_i, f_j})$, which will be unitary if both bases are normalized.  Thinking of $M$ with $1 \times (n-1)$ blocking,
\[
M =  \left[\begin{array}{c|c} m_{1,1} & r_1 \\ \hline c_1 & D\end{array}\right],
\]
the condition says that each element in the lower-right block $D$ has the same argument as the product of the first element in its row and the first element in its column.
\end{rmk}

\begin{rmk}
If a matrix is verified to be UECSM via Theorem \ref{test}, then an explicit conjugation $C$ is provided in terms of the eigenvectors of $A$ and $B$ via the proofs of Lemma \ref{realangles} and Theorem \ref{test}.  Let $\{e_i\}$ and $\{f_i\}$ be the orthonormal bases from Lemma \ref{realangles}, $U$ be the transition matrix from the standard basis to either $\{e_i\}$ or $\{f_i\}$, and let $J$ be the coordinate-wise complex conjugation operator on $\Cn$.  It follows from \cite[Sec. 3.2]{GarciaPutinar} that $U^*TU$ is a CSM and that
\[
C = UJU^* = UU^tJ.
\]
is a conjugation with respect to which $T$ is $C$-symmetric.
\end{rmk}

While there are already multiple proofs that all $2 \times 2$ matrices are UECSM (\cite[Ex. 6]{GarciaPutinar},\cite[Cor 3.3]{CFT},\cite[Cor. 1]{GWNewClasses}), we can use Theorem \ref{test} to provide another proof while explicitly exhibiting a conjugation $C$ for which $T=CT^*C$.

\begin{cor}
Every $2 \times 2$ matrix is unitarily equivalent to a complex symmetric matrix.
\end{cor}
\begin{proof}
Let $T = A + iB$ be a $2 \times 2$ matrix.  We may assume without loss of generality that $A$ and $B$ have 2 distinct eigenvalues, and that they do not share an eigenvector, as in either case $T$ would be normal and therefore UECSM by Corollary \ref{normal}.  Let $\{e_i\}$ and $\{f_i\}$ be orthonormal bases of eigenvectors of $A$ and $B$, respectively.  Compute the unitary $U = (u_{i,j}) = (\ip{e_i,f_j})$, and observe that since $A$ and $B$ do not share an eigenvector, no entry of $U$ is equal to $0$.  Multiply $e_1$ by a unimodular constant so that $\ip{e_1,f_1} \in \R$.  Since the columns of $U$ must be orthogonal, \mbox{$u_{1,1}\overline{u_{1,2}} + u_{2,1}\overline{u_{2,2}} = 0$} which yields
\[
\frac{u_{2,2}}{u_{1,2}u_{2,1}} = - \frac{u_{1,1}}{\lvert u_{1,2} \rvert^2}.
\]
As $u_{1,1} \in \R$, it follows that $\frac{u_{2,2}}{u_{1,2}u_{2,1}} \in \R$. By Theorem \ref{test}, $T$ is UECSM.  If we multiply $e_2$ by a unimodular constant so that $\ip{e_2, f_1} \in \R$ and let $V = \left[\begin{array}{c|c} e_1 & e_2\end{array}\right]$, then $T$ is $C$-symmetric for $Cx = VV^t\overline{x}$ by the remark following Theorem \ref{test}.
\end{proof}

\section{An algorithm for 3 $\times$ 3 matrices}
In this section we introduce an algorithm for determining whether or not a given $3 \times 3$ matrix is unitarily equivalent to a complex symmetric matrix.  We first require a few preparatory results.  The following proposition allows us to easily answer affirmatively in certain cases.

\begin{prop}\label{easycso}
Let $T = A + iB$ be a $3 \times 3$ matrix.  If either of the following conditions hold, then $T$ is unitarily equivalent to a complex symmetric matrix:
\begin{enumerate}\addtolength{\itemsep}{0.4\baselineskip}
\item[(i)] $A$ or $B$ has a repeated eigenvalue
\item[(ii)] $A$ and $B$ share an eigenvector
\end{enumerate}
\end{prop}
\begin{proof}
First (i).  Note that $T$ is UECSM if and only if $T - \lambda I$ is.  If $B$ has a repeated eigenvalue $\lambda$, then $B - \lambda I$ either has rank $0$ or $1$.  If $B - \lambda I = 0$, then $T - i\lambda I = A$, which is self-adjoint and therefore UECSM by Corollary \ref{normal}.  If $B - \lambda I$ has rank 1, we can show that $T$ is UECSM by simplifying an earlier result \mbox{\cite[Cor. 5]{GWNewClasses}} for the finite dimensional case. Let $v$ be a unit vector that spans $\Ran (B - \lambda I)$, and let $\M$ be the cyclic subspace generated by $v$ under $A$.  Clearly $\M$ is invariant under $A$, and it is easy to check that $\M^\perp$ is as well.  $T\restriction_{\M^\perp}$ is self-adjoint and $T\restriction_\M$ is cyclic self-adjoint, so by the Spectral Theorem we can assume without loss of generality that $v=(1,1,1)$ and that $A$ is diagonal \cite[Thm. 2.11.2]{Conway}.  Since $\Ran(B-\lambda I) = \Span \{v\}$, we have
\[
B = b\left[\begin{array}{ccc}1 & 1 & 1 \\ 1 & 1 & 1 \\ 1 & 1 & 1\end{array}\right]
\]
for some $b \in \R$.  This yields that $T=A + iB$ is UECSM.  The proof is similar if $A$ has a repeated eigenvalue.\\

Now (ii).  If $A$ and $B$ share an eigenvector, then up to unitary equivalence we know that $A = A_1 \oplus A_2$ and $B = B_1 \oplus B_2$ where $A_1$ and $B_1$ are $1 \times 1$ matrices and $A_2$ and $B_2$ are $2 \times 2$ matrices.  This yields that $T = (A_1 + iB_1) \oplus (A_2 + iB_2)$,  and because all $1 \times 1$ and $2 \times 2$ matrices are UECSM, $T$ is UECSM.
\end{proof}

The next lemma tells us that when Proposition \ref{easycso} does not imply that $T=A+iB$ is unitarily equivalent to a CSM, it is easy to construct a proper pair of orthogonal bases to which one can apply Theorem \ref{test}.

\begin{lemma}\label{makeproper}
If $T = A + iB$ is a $3 \times 3$ matrix which does not satisfy either hypothesis of Proposition \ref{easycso}, then any two orthogonal bases of eigenvectors $\{e_i\}$ and $\{f_i\}$ of $A$ and $B$, respectively, can be made proper by reordering them and scaling $e_1$.
\end{lemma}
\begin{proof}
Without loss of generality we may assume that $\norm{e_i} = \norm{f_i} = 1$, since the conditions of being proper are not affected by multiplying the basis vectors by real scalars.  Note that no $e_i$ is orthogonal to more than one $f_j$, since otherwise it would be a scalar multiple of the third element of $\{f_j\}$, a contradiction.  Similarly no fixed $f_i$ is orthogonal to more than one $e_i$.  In terms of the unitary matrix $U = (u_{i,j}) = (\ip{e_i,f_j})$, this means that no row or column has more than one $0$.\\

We claim that there is at most one $0$ in $U$.  If $U$ had more than one $0$, they must be in different rows and columns, so we could reorder the bases so that 
\[
U = \left[\begin{array}{ccc} 0 & * & * \\ * & 0 & * \\ a & b & *\end{array}\right],
\]
where the $*$'s represent arbitrary complex numbers.  To preserve the orthogonality of columns, we must have $a = 0$ or $b = 0$.  In either case, there must be more than one $0$ in a single column, which we have already excluded.  By reordering the bases, we can ensure that the $0$ entry is not in the first row or column.  If $\ip{e_1,f_1} \not \in \R$, multiply $e_1$ by $\frac{|\ip{e_1,f_1}|}{\ip{e_1,f_1}}$ and then $\{e_i\}$ and $\{f_i\}$ will be a proper pair.
\end{proof}

Using the preceding results, we can construct an algorithm that will decide whether or not a $3 \times 3$ matrix $T$ is unitarily equivalent to a complex symmetric matrix.  Since none of the operations are more complicated than finding roots of cubic polynomials, it can be performed using exact values, assuming the data is given exactly.  The author implemented it in \texttt{Mathematica} without much difficulty.  The algorithm is:
\begin{alg} Given a $3 \times 3$ matrix T,\\
\begin{enumerate}\addtolength{\itemsep}{0.5\baselineskip}
\item Compute $A =  \frac12(T + T^*)$ and $B = \frac1{2i}(T - T^*)$.
\item Compute the eigenvalues of $A$ and $B$.  If either $A$ or $B$ has a repeated eigenvalue, then $T$ is UECSM by Proposition \ref{easycso}.
\item Compute arbitrary sets of eigenvectors $\{g_i\}$ and $\{h_i\}$ of $A$ and $B$, respectively, and compute the matrix \mbox{$M = (m_{i,j}) = (\ip{g_i,h_j})_{i,j}$}.
\item If $M$ has more than one entry equal to $0$, then $T$ is UECSM (Lemma \ref{makeproper}, Proposition \ref{easycso}).  Otherwise, reorder the rows and columns of $M$ so that the $0$ entry is not in the first row or column (so that $\{g_i\}$ and $\{h_i\}$ form a proper pair).  Scale $g_1$ by $\frac{\abs{\ip{g_1,h_1}}}{\ip{g_1,h_1}}$.
\item By Theorem \ref{test}, $T$ is UECSM if and only if for all $2 \le i,j \le 3$,
\[\frac{m_{i,j}}{m_{1,j}m_{i,1}} \in \R.\]
\item If $T$ is UECSM, one can exhibit a corresponding conjugation $C$ by first normalizing $\{g_i\}$ and scaling $g_2$ and $g_3$ so that $\ip{g_i,h_1} \in \R$ for all $i$.  If \mbox{$U = \left[\begin{array}{c|c|c}g_1 & g_2 & g_3\end{array}\right]$}, then by the remark following Theorem \ref{test} $U^*TU$ is complex symmetric and $T$ is $C$-symmetric with respect to the conjugation \mbox{$Cx = UU^t\overline{x}$}.
\end{enumerate}
\end{alg}

It is worth noting that steps 1, 3 and 5 carry through to the $n \times n$ case as long as $A$ and $B$ have $n$ distinct eigenvalues and a proper pair of bases can be found.  Step 4 is no longer valid, as for $n > 3$ the preceding conditions do not guarantee that $T$ is UECSM.  A generalization of this algorithm to $n \times n$ matrices will be discussed in Section \ref{nbynsection}. The following examples illustrate the steps of the algorithm.

\begin{ex}\label{JordanWithParameter}
Theorem \ref{test} provides another proof of the fact from \cite[Ex. 1]{GWNewClasses} that for $a,b \ne 0$ the matrix
\[
T = \left[\begin{array}{ccc}0 & b & 0 \\ 0 & 0 & a \\ 0 & 0 & 0\end{array}\right]
\]
is unitarily equivalent to a complex symmetric matrix if and only if $\abs{a} = \abs{b}$.  By dividing by $b$, it is enough to consider matrices of the form
\[
T = \left[\begin{array}{ccc}0 & 1 & 0 \\ 0 & 0 & a \\ 0 & 0 & 0\end{array}\right].
\]  For this $T$, one can verify that
\[
\begin{array}{ccc}
A = \left[\begin{array}{ccc}0 & \frac12 & 0\\ \frac12 & 0 & \frac{a}2 \\ 0 & \frac{\overline{a}}2 & 0 \end{array}\right]
& \mbox{and} &
B = \left[\begin{array}{ccc}0 & -\frac{i}2 & 0\\ \frac{i}2 & 0 & -\frac{ia}2 \\ 0 & \frac{i\overline{a}}2 & 0 \end{array}\right].
\end{array}
\]
Furthermore, the eigenvalues of both $A$ and $B$ are $\{0, \frac12\sqrt{1 + \abs{a}^2},-\frac12\sqrt{1 + \abs{a}^2}\}$, the eigenvectors of $A$ are
\[\begin{array}{ccc}
g_1 = \left[\begin{array}{c}a \\ 0 \\ -1\end{array}\right], &
g_2 = \left[\begin{array}{c}1 \\ \sqrt{1+\abs{a}^2} \\ \overline{a}\end{array}\right], &
g_3 = \left[\begin{array}{c}1 \\ -\sqrt{1+\abs{a}^2} \\ \overline{a}\end{array}\right],
\end{array}\]
and the eigenvectors of $B$ are
\[\begin{array}{ccc}
h_1 = \left[\begin{array}{c}a \\ 0 \\ 1\end{array}\right], &
h_2 = \left[\begin{array}{c}1 \\ i\sqrt{1+\abs{a}^2} \\ -\overline{a}\end{array}\right], &
h_3 = \left[\begin{array}{c}1 \\ -i\sqrt{1+\abs{a}^2} \\ -\overline{a}\end{array}\right].
\end{array}\]
The matrix $M=(\ip{g_i,h_j})_{i,j=1}^3$ required by the algorithm is given by:
\[
M = \left[\begin{array}{ccc}1-\abs{a}^2 & 2 & 2 \\ 2 & \beta & \overline{\beta} \\ 2 & \overline{\beta} & \beta\end{array}\right].
\]
where $\beta = 1+i - \frac{1-i}{\abs{a}^2}$.  If $\abs{a} \ne 1$, then the bases $\{g_i\}$ and $\{h_i\}$ are proper and since $\beta$ has non-zero imaginary component, $T$ is not UECSM by Theorem \ref{test}.\\

If $\abs{a} = 1$, then $\beta = 2i$.  We can relabel the vectors of each basis and scale the new $e_1$ by $-i$ to get the matrix,
\[
M^\prime = 2\left[\begin{array}{ccc}1 & -1 & -i \\ -i & i & 1 \\ 1 & 1 & 0\end{array}\right]
\]
It is easy to check that $\frac{m_{i,j}}{m_{i,1}m_{1,j}} \in \R$ for $i,j \ge 2$, so $T$ is UECSM by Theorem \ref{test}.
\end{ex}

\begin{ex}
Let
\[
T=
\left[
\begin{array}{ccc}
 1+4 i & (-2-i) \sqrt{2} & -1-4 i \\
 i \sqrt{2} & 0 & i \sqrt{2} \\
 -1 & (2-i) \sqrt{2} & 1
\end{array}
\right].
\]

In this example, we prove that $T=A+iB$ is unitarily equivalent to a complex symmetric matrix and use the full algorithm to find a conjugation $C$ with respect to which $T$ is $C$-symmetric.  Per step 2, we first calculate the eigenvalues of $A$ and $B$, which are 
\[
\left\{2 \left(1+\sqrt{2}\right),-2,2 \left(1-\sqrt{2}\right)\right\}\quad \mbox{and} \quad  
\left\{2 \left(1+\sqrt{3}\right),2 \left(1-\sqrt{3}\right),0\right\},
\]
respectively.  Neither $A$ nor $B$ has a repeated eigenvalue, and so we calculate the eigenvectors of $A$:
\[\begin{array}{ccc}
g_1 = \left[\begin{array}{c}-1-2 i \sqrt{2} \\ 2+i \sqrt{2} \\ 3 \end{array} \right],&
g_2 = \left[\begin{array}{c}1 \\ -i \sqrt{2} \\ 1\end{array}\right],&
g_3 = \left[\begin{array}{c}-1+2 i \sqrt{2} \\ -2+i\sqrt{2} \\ 3 \end{array} \right],
\end{array}\]
and the eigenvectors of $B$:
\[\begin{array}{ccc}
h_1=\left[\begin{array}{c}-1-\frac{2}{\sqrt{3}} \\i\sqrt{\frac{2}{3}} \\ 1\end{array}\right],&
h_2=\left[\begin{array}{c}-1+\frac{2}{\sqrt{3}} \\-i\sqrt{\frac{2}{3}} \\ 1\end{array}\right],&
h_3=\left[\begin{array}{c}1 \\i \sqrt{2} \\1\end{array}\right].
\end{array}\]

Since $A$ and $B$ do not share an eigenvector, we must use the full algorithm and not a shortcut provided by Proposition \ref{easycso}.  Next we compute $M=(\ip{g_i,h_j})$:
\[
M =
\left[
\begin{array}{ccc}
 \frac{2}{3} \left(2+i \sqrt{2}\right) \left(3+\sqrt{3}\right) & -\frac{2}{3} i \left(-2
   i+\sqrt{2}\right) \left(-3+\sqrt{3}\right) & 4-4 i \sqrt{2} \\
 -\frac{4}{\sqrt{3}} & \frac{4}{\sqrt{3}} & 0 \\
 \frac{2}{3} \left(2-i \sqrt{2}\right) \left(3+\sqrt{3}\right) & \frac{2}{3} i \left(2
   i+\sqrt{2}\right) \left(-3+\sqrt{3}\right) & 4+4 i \sqrt{2}
\end{array}
\right].
\]
If we let $\alpha = \frac{\abs{m_{1,1}}}{m_{1,1}}$ and $\beta = \frac{\abs{m_{3,1}}}{m_{3,1}}$, rewriting $M$ with respect to $\{\alpha g_1, g_2, \beta g_3\}$ and $\{h_1, h_2, -ih_3\}$ we get,
\[
M^\prime = \left[
\begin{array}{ccc}
 2 \left(\sqrt{2}+\sqrt{6}\right) & 2 \sqrt{2} \left(-1+\sqrt{3}\right) & 4 \sqrt{3} \\
 -\frac{4}{\sqrt{3}} & \frac{4}{\sqrt{3}} & 0 \\
 2 \left(\sqrt{2}+\sqrt{6}\right) & 2 \sqrt{2} \left(-1+\sqrt{3}\right) & -4 \sqrt{3}
\end{array}
\right].
\]

As all of the entries are real, Theorem $1$ says that $T$ is UECSM.  Letting $\{e_i\} = \{\frac{\alpha g_1}{\norm{\alpha g_1}}, \frac{g_2}{\norm{g_2}}, \frac{\beta g_3}{\norm{\beta g_3}}\}$ and $U = \left[e_1 | e_2 | e_3 \right]$, we can construct a conjugation $C$ for which $T$ is $C$-symmetric and a complex symmetric matrix that $T$ is unitarily equivalent to:
\[
Cx = UU^t\overline{x} = 
\left[\begin{array}{ccc}
 \frac{1}{2} & -\frac{i}{\sqrt{2}} & -\frac{1}{2} \\
 -\frac{i}{\sqrt{2}} & 0 & -\frac{i}{\sqrt{2}} \\
 -\frac{1}{2} & -\frac{i}{\sqrt{2}} & \frac{1}{2}
\end{array}\right]
\left[\begin{array}{c}\overline{x_1} \\ \overline{x_2} \\ \overline{x_3}\end{array}\right],
\]
\[
U^*TU = 2\left[
\begin{array}{ccc}
  1+\sqrt{2}+i  & -i & i \\
 -i & -1 & -i \\
 i & -i & 1-\sqrt{2}+i
\end{array}
\right].
\]

\end{ex}

\begin{ex}\label{intex}
Let\[
T_1 = \left[\begin{array}{ccc}0 & 7 & 0 \\ 0 & 1 & -5 \\ 0 & 0 & 6\end{array}\right], \quad
T_2 = \left[\begin{array}{ccc}0 & 7 & 0 \\ 0 & 1 & -5 \\ 0 & 0 & 3\end{array}\right].
\]
Applying the same method as in the previous example yields that $T_1$ is $C$-symmetric where
\[
C\left[\begin{array}{c}\vspace{0.03in}x_1 \\ \vspace{0.03in}x_2 \\ \vspace{0.03in}x_3\end{array}\right] = 
\left[
\begin{array}{ccc}
  \vspace{0.03in}\frac{6 \left(-19+ 6i\sqrt{74}\right)}{3025} &
   \frac{42 \left(19-6 i \sqrt{74}\right)}{3025} &
   \frac{7}{605} \left(19-6 i \sqrt{74}\right) \\
 \vspace{0.03in}\frac{42 \left(19-6 i \sqrt{74}\right)}{3025} &
   \frac{19 \left(-19+6i\sqrt{74}\right)}{3025} &
   \frac{6}{605} \left(19-6 i \sqrt{74}\right) \\
 \vspace{0.03in}\frac{7}{605} \left(19-6 i \sqrt{74}\right) &
   \frac{6}{605} \left(19-6 i \sqrt{74}\right) &
   \frac{6}{605} \left(-19+6i\sqrt{74}\right)
\end{array}
\right]
\left[\begin{array}{c}\vspace{0.03in}\overline{x_1} \\ \vspace{0.03in}\overline{x_2} \\ \vspace{0.03in}\overline{x_3}\end{array}\right],
\]
and that $T_1$ is unitarily equivalent to
\[
T_1^\prime = 
\left[
\begin{array}{ccc}
 \frac{56}{37}-i \sqrt{\frac{37}{2}} & -\frac{55}{37} & \frac{35 \sqrt{55}}{74} \\
 -\frac{55}{37} & \frac{56}{37}+i \sqrt{\frac{37}{2}} & \frac{35 \sqrt{55}}{74} \\
 \frac{35 \sqrt{55}}{74} & \frac{35 \sqrt{55}}{74} & \frac{147}{37}
\end{array}
\right].
\]

We also get that $T_2$ is not UECSM.  It is easy to check that if $T_2$ has Cartesian decomposition $T_2 = A + iB$, both $A$ and $B$ have $3$ distinct eigenvalues and that they do not share an eigenvector.  None of the eigenspaces of $A$ are orthogonal to any of the eigenspaces of $B$, and so it is easy to construct a proper pair of orthogonal bases of eigenvectors of $A$ and $B$.  It is also easy to check that these will not satisfy condition \eqref{realratio} of Theorem \ref{test}.
\end{ex}

\section{Applications for generic $n \times n$ matrices}\label{nbynsection}
Theorem \ref{test} provides a sufficient condition for an arbitrary $n \times n$ matrix to be unitarily equivalent to a complex symmetric matrix, so long as a proper pair of bases can be found, and if $A$ and $B$ each have $n$ distinct eigenvalues then Theorem \ref{test} is necessary as well as sufficient.  It is not guaranteed that an arbitrary pair of orthogonal bases can be made proper simply by reordering and scaling the elements.  However, for a random matrix $T$, its Cartesian components $A$ and $B$ will both have $n$ distinct eigenvalues with probability one \cite[Chapter 3]{Mehta}.  Moreover, with probability one an arbitrary choice of bases of eigenvectors of $A$ and $B$ can easily be made proper, as in Lemma \ref{makeproper}.  For more on random self-adjoint matrices, one can consult \cite{Mehta}.  In light of the above, the algorithm given in the previous section will almost surely apply to a random $n \times n$ matrix, which is useful for probabilistic searches.\\

For instance, in \cite[Sec. 5]{GWNewClasses} the authors determine that any $4 \times 4$ partial isometry $T$ with $\rank T = 1$ or $\rank T = 3$ is UECSM.  $T$ is also trivially UECSM in the cases where $\rank T = 0$ or $\rank T = 4$.  However, they were unable to answer whether or not every $4 \times 4$ partial isometry with $\rank T = 2$ is unitarily equivalent to a complex symmetric matrix.  All such partial isometries are unitarily equivalent to a matrix of the form $UP$, where $P$ the projection onto the first two standard basis vectors and $U$ is unitary.  We used \texttt{Mathematica} to test $100,000$ matrices of the form $UP$, with random unitary components generated by taking the matrix exponential of a random skew-Hermitian matrix.  All of them were UECSM, and so we conjecture:
\begin{conjecture}
Every rank-two $4 \times 4$ partial isometry is unitarily equivalent to a complex symmetric matrix.
\end{conjecture}

\section*{Acknowledgements}
We wish to thank S.R. Garcia for his many valuable suggestions, continued support, and financial contribution via NSF Grant DMS 0638789.

\end{document}